\documentclass[11pt]{article}
\usepackage{amssymb}
\usepackage{enumerate}
\usepackage{enumitem}
\usepackage{mathrsfs}
\usepackage{latexsym,amsmath,amssymb,amsfonts,epsfig,graphicx,cite,psfrag}
\usepackage{eepic,color,colordvi,amscd}
\usepackage{verbatim}
\usepackage{appendix}

\usepackage{float}
\usepackage{graphicx}
\usepackage{cite}
\usepackage{caption}
\usepackage{subcaption}

\newcommand{\qed}{\hfill $\Box $}

\newtheorem{theorem}{Theorem}
\newtheorem{lemma}[theorem]{Lemma}

\newcommand{\ar}{{\rm ar}}
\newcommand{\ex}{{\rm ex}}
\newcommand{\G}{{\cal G}}

\begin{document}
	
	\title{Anti-Ramsey number of intersecting cliques\footnote{Research supported by National Key Research and Development Program of China 2023YFA1010200 and  National Natural
Science Foundation of China grant  12271425.}}
	
	\author{
Hongliang Lu, Xinyue Luo and Xinxin Ma \\School of Mathematics and Statistics\\
Xi'an Jiaotong University\\
Xi'an, Shaanxi 710049, China
}
\date{}

%\underline{}\maketitle

\date{}

\maketitle
	\begin{abstract}

An edge-colored graph is called a rainbow graph if all its edges have distinct colors. The anti-Ramsey number \( ar(n, G) \), for a graph \( G \) and a positive integer \( n \), is defined as the minimum number of colors \( r \) such that every exact \( r \)-edge-coloring of the complete graph \( K_n \) contains at least one rainbow copy of \( G \). A \( (k, r) \)-fan graph, denoted \( F_{k, r} \), is a graph composed of \( k \) cliques each of size \( r \), all intersecting at exactly one common vertex.
In this paper, we determine \( ar(n, F_{k, r}) \) for  \( n \geq 256r^{16}(k+1)^5 \), \( k \geq 1 \), and \( r \geq 2 \). %Our result generalizes previous work on anti-Ramsey numbers for \( ar(n, F_{k, 3}) \)
		
	\end{abstract}
	\begin{flushleft}
		%{\em AMS classification:} 05C70
	\end{flushleft}
	
	\section{Introduction}
	
All graphs considered in this paper are finite and simple, meaning they have no multiple edges or loops. We use \( G=(V(G), E(G)) \) to denote a graph with vertex set $V(G)$  and edge set $E(G)$. The order of \( G \), denoted by \( |G| \), is defined as the cardinality of the vertex set \( V(G) \), that is \( |G| := |V(G)| \), and the size of \( G \), denoted by \( e(G) \), is the cardinality of the edge set \( E(G) \), so \( e(G) := |E(G)| \). For any subset \( X\subseteq V(G) \), we use \( G - X \) to represent the subgraph of \( G \) obtained by deleting all vertices in \( X \) and all edges incident to any vertex in \( X \). When \( X=\{x\} \), we simply write it as \( G - x \). For \( Y\subseteq E(G) \), \( G - Y \) represents the subgraph obtained from \( G \) by removing all edges in \( Y \). For an edge set \( Y \) such that \( Y\cap E(G)=\emptyset \), we let \( G + Y \) denote the graph whose vertex set is \( V(G)\cup V(Y) \) and whose edge set is \( E(G)\cup E(Y) \).
For two nonempty subsets \( S,T\subseteq V(G) \), let \( E_G(S,T) \) denote the set of all edges in \( G \) with one endpoint in \( S \) and one endpoint in \( T \).
The neighborhood of vertex \( x \) in \( G \) is denoted by \( N_G(x) \). The \emph{degree} of \( x \) in \( G \), denoted by \( d_G(x) \), is equal to the size of \( N_G(x) \).
We use $\delta(G)$ and $ \Delta(G) $ to denote the minimum and maximum degrees of the vertices of $ G $. A \emph{matching} in \( G \) is a set of edges from \( E(G) \) such that no two edges share a common vertex. The \textit{matching number} of $ G $, denoted by $ \nu(G) $, is the maximum number of edges in a matching in $ G $.

Given two graphs \( H \) and \( G \), if \( G \) does not contain \( H \) as  a subgraph, we say that \( G \) is \( H \)\emph{-free}.
For integers \( n \geq p \geq 1 \), let \( T_{n,p} \) denote the Tur\'an graph, i.e., the complete \( p \)-partite graph on \( n \) vertices where each partition class has either \( \lfloor n/p \rfloor \) or \( \lceil n/p \rceil \) vertices, with edges connecting all pairs of vertices from different parts. Let \( t_p(n) \) denote the number of edges in \( T_{n,p} \). Let $K_n$ be the complete graph with $ n $ vertices. Given two positive integers $r,k$, let $kK_r$ denote the union of $k$ vertex-disjoint complete graph $K_r$. Let $[n]:=\{1,2, \cdots ,n \}$.
For two vertex-disjoint graphs $G$ and $H$, the \textit{join} of $G$ and $H$, denoted by $G \lor H$, is the graph obtained from $G \cup H$ by adding edges joining every vertex of $G$ to every vertex of $H$. The \textit{union} of $G$ and $H$ is the graph $G \cup H$ with vertex set $V(G) \cup V(H)$ and edge set $E(G) \cup E(H)$.

	An $r$-\textit{edge-coloring} of a graph is an assignment of $r$ colors to the edges of the graph. An \textit{exactly $r$-edge-coloring} of a graph is an $r$-edge-coloring that uses all $r$ colors. An edge-colored graph is called \textit{rainbow} if all edges have distinct colors. Given a positive integer $ n $ and a graph set $ \mathcal{G} $ , the \textit{anti-Ramsey number} $ ar(n,\mathcal{G}) $ is the minimum number of colors $r$ such that each edge-coloring of $K_n$ with exactly $r$ colors contains a rainbow copy of $G \in \mathcal{G}$.
	%there exists a rainbow copy of $G$ for some $G \in \mathcal{G}$ in any exactly $r$-edge-coloring of $K_n$.
	When $\mathcal{G}=\{G\}$, we denote $ar(n,\mathcal{G})$  by $ ar(n,G) $ for simplicity. The \textit{Tur\'an number} $ \ex(n,G) $ is the maximum number of edges of a graph on $ n $ vertices containing no subgraph isomorphic to $ G $. 	 A
	graph $H$ on $n$ vertices with $\ex(n, G)$ edges and without a copy of $G$ is called  \emph{extremal
	graph} for $G$-free graphs. We use $EX(n,G)$ to denote the set of extremal graphs for $G$, i.e.
	\[
	EX(n,G):=\{H\ |\ |V(H)|=n, e(H)= \ex(n,G), \mbox{ and $H$ is $G$-free}\}.
	\]

	The value of $ar(n,G)$ is closely related to the Tur\'an number $\ex(n,G)$ as the following inequality
	\begin{align}\label{low-upp}
		2+\ex(n,\G) \leq ar(n,G) \leq 1+\ex(n,G),
	\end{align}
	where $\G=\{G-e\ |\ e\in E(G)\}$.  Erd\H{o}s, Simonovits and S\'os \cite{Erdos1973}  proved there exists a number $n_0(p)$ such that $ar(n,K_p)=t_{p-1}(n)+2$ for $n> n_0(p)$. Later, Montellano-Ballesteros \cite{Mon-Ball2005} and, independently, Neumann-Lara \cite{Mon-Ball2002} extended this result to all $n$ and $p$ satisfying $n> p\geq 3$.
	
Another interesting problem in anti-Ramsey number theory is determining the anti-Ramsey numbers for cycles. Let $l\geq 3$ be a positive integer and let $C_l$ denote the cycle of length $l$.  Erd\H{o}s, Simonovits and S\'os \cite{Erdos1973} conjectured that
\[
ar(n,C_l) = \left( \frac{l-2}{2} + \frac{1}{l-1} \right) n + O\left( 1 \right)
\]
and verified the conjecture for $l=3$. Alon \cite{Alon1983} confirm the case for $l=4$ by showing that $ar(n,C_4)=\lfloor \frac{4n}{3} \rfloor-1$. Jiang, Schiermeyer and West \cite{Jiang} extended the proof to all $l \leq 7$. The conjecture was fully resolved by Montellano-Ballesteros and Neumann-Lara \cite{Mon-Ball2005}.

	The Tur\'an number  for matchings is  determined by Erd\H{o}s and Gallai\cite{erdos1959} for $ n\geq 2t\geq 2 $.  For anti-Ramsey numbers of matchings, Schiermeyer\cite{Schiermeyer2004} first established that  $ \ar(n,tK_2)=\ex(n,(t-1)K_2)+2 $  for $ n\geq 3t+3$.   This result was later improved by Fujita et al. \cite{Fujita}, who extended the range to \( n \geq 2s + 1 \). A complete resolution for all \( s \geq 2 \) and \( n \geq 2s \) was achieved by Chen, Li, and Tu \cite{Chen2009}, who determined the exact value of \( ar(n, M_s) \). Interestingly,  Haas and Young \cite{Haas2012} gave a simple proof for the case $n=2s$. There is a large volume of
	literature on the anti-Ramsey number of graphs. Interested readers are referred to the survey
	by Fujita, Magnant, and Ozeki \cite{Fujita}.

%Later, Chen, Li and Tu \cite{Chen2009} and independently Fujita, Kaneko, Schiermeyer and Suzuki \cite{Fujita} showed that $ \ar(n,tK_2)=\ex(n,(t-1)K_2)+2 $ for $ n\geq 2t+1\geq 5 $. The values $ \ar(2t,tK_2)=\ex(2t,(t-1)K_2)+2 $ for $ 2\leq t\leq 6 $ and  $ \ar(2t,tK_2)=\ex(2t,(t-1)K_2)+3 $ for $ t\geq 7 $ were determined in and by Haas and Young \cite{Haas2012}, independently. %Gilbos and Roditty \cite{Gilboa2016} considered the graph with small connected components and proved inductive results of the form ``if $ \ar(n,G\cup t_0P_s)\leq %f(n,t_0,G) $ for sufficiently large $ n $ where $ p_s $ is a path of length $ s-1 $, then $ \ar(n,G\cup tP_s)\leq f(n,t,G) $ for sufficiently large $ n $ and $ t\geq t_0 $, where $ s=2 $ or $ 3 $".

 %Yuan and Zhang \cite{Yuan} provided the exact results of $ar(n,kC_3)$ when $n$ is sufficiently large. Later, Wu and Zhang et al. \cite{Wu} improved the result from $n$ is sufficiently large to $n \geq 2k^2-k+2$.	Lu, Luo and Ma \cite{} obtained the new bound for $n\geq 15k+57$.
	
	A graph on $ (r-1)k+1 $ vertices consisting of $ k $ cliques each with $ r $ vertices, which intersect in exactly one common vertex, is called a \textit{$(k,r) $-fan}  and denoted by $ F_{k,r} $. Liu, Lu and Luo \cite{LLL} determined the anti-Ramsey number $\ar(n,F_{k,3})$ for $n\geq 100k^2$. In this paper, we determined  $\ar(n,F_{k,r})$ for large $n$.
	
	\begin{theorem}\label{main}
		Let $ n,k,r $ be integers such that $ k\geq 1 $ and $ r\geq 2$,  and  $ n\geq 256(k+1)^5r^{16} $. The following holds $$ ar(n,F_{k+1,r})=\ex(n,F_{k,r})+2.$$
	\end{theorem}
	
	The rest of the paper is organized as follows. In Section 2, we introduce some technical lemmas. In Section 3, we complete the proof of Lemma \ref{partition}. In Section 4, we give the proof of Theorem \ref{main}.
	
	\section{Several Technical Lemmas}
	For two integers $\nu, \Delta \geq 1$, define $ f\left(\nu, \Delta\right)=\max\{e(G)\ |\ \nu(G)\leq \nu, \Delta(G)\leq \Delta\} $. Let $\mathcal{F}_{\nu, \Delta}:=\{G\ |\ \nu(G)\leq \nu, \Delta(G)\leq \Delta, e(G)=f\left(\nu, \Delta\right)\}$. In 1972, Abbott, Hanson and Sauer\cite{Abbott} determined $ f(k, k)$.
	\begin{theorem}[Abbott, Hanson and Sauer, \cite{Abbott}]\label{kk-value}
		Let $k\geq 2$ be an integer. Then
		\begin{equation*}
			f(k,k)=\begin{cases}
				k^2+(k-1)/2 ,& \mbox{if k is odd}, \medskip\\
				k^2+k,& \mbox{if k is even}.
			\end{cases}	
		\end{equation*}
	\end{theorem}

	Later Chv\'atal and Hanson\cite{chva1976}  proved the following result, which generalized Theorem \ref{kk-value}.
	\begin{theorem}[Chv\'atal and Hanson, \cite{chva1976}]\label{Chvatal}
		For every $\nu \geq 1$ and $\Delta\geq 1$,
		$$
		f( \nu ,\Delta ) =\nu \Delta +\lfloor \frac{\Delta}{2} \rfloor \lfloor \frac{\nu}{\lceil \Delta/2 \rceil} \rfloor \leq \nu \Delta +\nu.
		$$
	\end{theorem}
	
	%
%	Since $ \ex(n,K_r)-\ex(n-1,K_r)=\lfloor (\frac{r-2}{r-1})n\rfloor$, the following lemma holds.
	\begin{lemma}[Chen et al., \cite{chen2003}]\label{delta}
		Let $ G $ be a graph of order $ n $. Let $ k $ be an integer and $ c $ some constant independent from $ n $. If $ e(G)\geq \ex(n,K_r)+c $ and $ d_G(x)\leq (\frac{r-2}{r-1})n-(k+1) $, then $ e(G-x)\geq \ex(n-1,K_r)+c+(k+1) $.
	\end{lemma}
	%A graph on $ 2k+1 $ vertices consisting of $ k $ triangles which intersect in exactly one	common vertex is called a $ k $-fan and denoted by $ F_k $, i.e. $ F_{k,3} $: For each $ k $,

  Erd\H{o}s et al. \cite{erdos1995} determined the $ \ex(n,F_{k,3})$ for large $n$.
	\begin{theorem}[Erd\H{o}s, F{\"u}redi, Gould and Gunderson, \cite{erdos1995}]\label{Fkextremal}
		For two integers $n,k$ such that $k\geq 2$ and $n\geq 50k^2 $,   $$
		\ex(n,F_k)=\left\{ \begin{matrix}
			\lfloor \frac{n^2}{4} \rfloor +k^2-k, &		\mbox{if $k$ is odd}, \medskip\\
			\lfloor \frac{n^2}{4} \rfloor +k^2-\frac{3}{2}k, &\mbox{if $k$ is even}.\\
		\end{matrix} \right.
		$$
		%edges, then $G$ contains a copy of $F_k$, %Furthermore, the number of edges is the best possible.
	\end{theorem}
	Erd\H{o}s, F{\"u}redi, Gould and Gunderson \cite{erdos1995} showed that each member in $ EX(n,F_{k,3}) $ could be obtained from a Tur\'an graph $ T_{n,2} $ with a graph $ G $ embedding in one partite set, where $ G $ is empty for $ k = 1 $ and $\mathcal{F}_{k-1, k-1} $ for $ k \geq 2 $.
	Chen et al. determined the $ \ex(n,F_{k,r})$ for large $n$.
	
	\begin{theorem}[Chen et al., \cite{chen2003}]\label{F_{k,r}extremal}
		For every $ k\geq 1 $ and $ r\geq 3 $, and for every $n\geq 16k^3r^8 $,   $$
		\ex(n,F_{k,r})=\left\{ \begin{matrix}
			\ex(n,K_r)  +k^2-k, &		\mbox{if $k$ is odd}, \medskip\\
			\ex(n,K_r)  +k^2-\frac{3}{2}k, &\mbox{if $k$ is even}.\\
		\end{matrix} \right.
		$$
	Each member in $ EX(n,F_{k,r}) $ may be obtained from a Tur\'an graph $ T_{n,r-1} $ with a graph $ G $ embedding in one partite set, where $ G $ is empty for $ k = 1 $ and $\mathcal{F}_{k-1, k-1} $ for $ k \geq 2 $.	%edges, then $G$ contains a copy of $F_k$, %Furthermore, the number of edges is the best possible.
	\end{theorem}
	%Each member in $ EX(n,F_{k,r}) $ may be obtained from a Tur\'an graph $ T_{n,r-1} $ with a graph $ G $ embedding in one partite set, where $ G $ is empty for $ k = 1 $ and $\mathcal{F}_{k-1, k-1} $ for $ k \geq 2 $.

	Let $ b $ be a positive integer and let $ X $ and $ Y $ be two disjoint vertex subsets of $ V(G) $: We say that $ X $ \textit{dominates} $Y$ with \textit{b-deficiency} if $d_Y(x)\geq |Y|-b $ for each $ x\in X $. Let $ V_1,..., V_m  $ be disjoint subsets of $ V(G) $. We say that $ \{V_1,..., V_m\} $ is \textit{b-deficiency complete} if $ V_i $ dominates $ V_j $ with deficiency $ b $ for every pair $ i\ne j $ with $ i,j=1,2,...,m $.
	
	\begin{lemma}[Chen et al., \cite{chen2003}]\label{tcycles}
		Let $ b $ be a positive integer. Let $G$ be a graph and let $\{X_1,...,X_m\}$ be a $b$-deficiency complete partition of $ V(G) $ with $ |X_i|\geq mb+2t $  for each $i$. Suppose that $ C_1,C_2,...,C_t $ are $ t $ cliques of $G$ with the properities:
		\item[(a)] $ |C_i\cap X_j|\leq 2 $ for each pair $ i $ and $ j $,
		\item[(b)] $ |C_i\cap X_j|= 2 $ at most one $ j $ for each $ i $.
		
		Then, there exist $ t $ cliques $ D_1,..,D_t $ satisfying:
		\item[(1)] $ C_i\subseteq D_i $ for each $ i $,
		\item[(2)] $ D_1-C_1,...,D_t-C_t $ are mutually disjoint,
		\item[(3)] for each $ i $ we have that $ |D_i\cap X_j|= 1 $ for all $j$ except possibly one at which $ |D_i\cap X_j|=|C_i\cap X_j|= 2 $.
	\end{lemma}
	
	\begin{lemma}[Chen et al., \cite{chen2003}]\label{k cycles}
		Let $G$ be a graph and $ V_1,V_2,...,V_m $ be $m$ vertex disjoint subsets of $ V(G) $ and $ Y_0\subseteq V(G)-\cup_{i=1}^mY_i $ such that $ |Y_i|\geq (i-1)b+k $ for each $ i=1,...,m $. If $ Y_j $ dominates $ Y_i $ with $b$-deficiency for each $ i=1,...,m, j=0,1,...,m $ and $i\ne j$, then, there are $k$ vertex disjoint cliques $ C_1,C_2,...,C_k $ satisfying $ |C_i|=m $ and $ |C_i\cap Y_j|=1 $ for each $i$ and $j\geq 1$. Furthermore, if $ |Y_0|\geq mb+k $, then there are $k$ vertex disjoint cliques $ D_1,D_2,...,D_k $ with the property that $ |D_i|=m+1 $ and $ |D_i\cap Y_j|=1 $ for each $i=1,..,k$ and $j=0,1,...,m$.
	\end{lemma}
	
	%We will heavily use the tools provided by Chen et al.\cite{chen2003} to prove Theorem \ref{main}.
For complete the proof of Theorem \ref{main}, we need the following lemma, whose proof employs some idea presented in \cite{chen2003}. %We  will be provided the proof  in Section 3.
	\begin{lemma}\label{partition}
		For every $ k\geq 1 $ and $ r\geq 4 $, suppose that $G$ is an $ F_{k+1,r} $-free graph on $n$ vertices % and $ ex(n,F_{k,r})+2 $ edges
		with $ n\geq 4(k+1)^2r^4  $ and with minimum degree $ \delta(G)> (\frac{r-2}{r-1})n-(k+1) $, and $ e(G)> ex(n,F_{k,r})$. Then there exists a partition $ V(G)=V_0\cup V_1\cup\cdots\cup V_{r-2} $, such that  for all $ i=0,...,r-2 $, the following three statements  hold:
\begin{itemize}
  \item [$(i)$]  $|V_i|\geq \frac{n}{r-1}-(k+1)(2r-3)$;
  \item [$(ii)$] %\begin{equation}\label{(1)}
		$	\sum_{j\ne i}\nu(G[V_j])\leq k\ \ and\ \  \Delta(G[V_i])\leq k;$
		%\end{equation}
  \item [$(iii)$] %\begin{equation}\label{(2)}
			$ d_{G[V_i]}(x)+\sum_{j\ne i}\nu(G[N(x)\cap V_j])\leq k\quad \mbox{for any $x\in V_i$}.$
		%\end{equation}
\end{itemize}
	\end{lemma}
	
	Let $G$ be a graph with a partition of the vertices into $ r-1 $ non-empty parts $V(G)=V_0\cup V_1\cup\cdots\cup V_{r-2}$. Define $G_i=G[V_i]$ for each $ i=0,1,...,r-2 $, and define
	$$
	G_{cr}=G-(\cup_{i=0}^{r-2}E(G[V_i])).%(V(G),\{ v_i v_j \in E(G)\ |\  v_i \in V_i, v_j \in V_j, i\ne j \}).
	$$
	\begin{lemma}[Chen et al., \cite{chen2003}]\label{edgelow}
		Suppose $ G $ is partitioned as above so that (\ref{(1)}) and (\ref{(2)}) are satisfied.  If $ G $ is $ F_{k+1,r} $-free, then
		\begin{equation*}
			\sum_{i=0}^{r-2}e(G_i)-(\sum_{0\leq i<j\leq r-2}|V_i||V_j|-e(G_{cr}))\leq f(k,k).
		\end{equation*}
	\end{lemma}
	
	\section{Proof of  Lemma \ref{partition}}
	\noindent \textbf{Proof of Lemma \ref{partition}.} Since $ e(G)>ex(n,F_{k,r}) $, then $ G $ contains $ k $ edge disjoint cliques $ D_1,...,D_k $ sharing one vertex $ v_0 $ with $ |D_i|=r $ for $ i=1,...,k $. Let $K$ be a clique of maximum size in $ G-(\cup_{i=1}^{k}V(D_i)\backslash\{v_0\}) $ containing the vertex $ v_0 $.
	%The clique $ K $ is chosen to be as large as possible under these constraints.
Denote $ V(K)=\{v_0,v_1,...,v_s\} $. From the inequality $ \delta(G)> (\frac{r-2}{r-1})n-(k+1) $, it follows that $ |K|\geq 2 $. So we have $ s\geq 1 $ and $ |\cup_{i=1}^kD_i\cup K|=k(r-1)+s+1 $. For each $ i=0,...,s $, let $ X_i=\cap_{j\ne i}N_G(v_j)-V(\cup_{j=1}^kD_j\cup K)$. Since $ K $ is a clique that contains $ v_0 $  and has the maximum order and is  vertex-disjoint with $ D_1,...,D_k $, it follows $$ X_i\cap X_j=\emptyset\   \mbox{for $ i\ne j $}.$$
	Note  that $ \delta(G)> (\frac{r-2}{r-1})n-(k+1) $. So we may infer that
\begin{align*}
	 |X_i\cup V(\cup_{z=1}^kD_z\cup K)|\geq n-\left( \frac{n}{r-1}+(k+1)\right)s,
\end{align*}
which implies that
\begin{align}\label{Xi-lower}
 |X_i|\geq\frac{r-1-s}{r-1}n-(k+1)s-k(r-1)-s-1. 
 \end{align}
	
\medskip
	\textbf{Claim 1.}~$X_i$ is an independent set for $i=1,\ldots,s$.
	\medskip

	Otherwise, suppose that there is an edge $ uv\in E(G[X_i]) $. Replacing $ v_i $ by the edge $ uv $ in $ D $, we obtain a copy of clique whose size greater than $K$, a contradiction. This completes the proof of Claim 1.
	
	\medskip
	\textbf{Claim 2.}~$ s=r-2 $.
	\medskip
	
By contradiction, suppose  that $s\neq r-2 $.	Since $ G $ is $ F_{k+1,r} $-free, we may assume that $s\leq r-3 $.
For $i\in [s]$, by Claim 1 and (\ref{Xi-lower}),
	\begin{align*}
		d_G(x_i)&\leq n-|X_i|\\
		&\leq n-\left( \frac{2}{r-1}n-(k+1)(r-3)-k(r-1)-(r-2)\right) \\
		&=\frac{r-3}{r-1}n+(k+1)(r-3)+k(r-1)+r-2.
	\end{align*}
	contradicts with the fact that $ \delta(G)>\frac{r-2}{r-1}n-(k+1)$.
	This completes the proof of Claim 2.
	
By Claim 2 and  (\ref{Xi-lower}), we may infer that
for each $ i=0,...,r-2 $,
	\begin{equation}\label{Xi-lower2}
		|X_i|\geq\frac{1}{r-1}n-(k+1)(2r-3).
	\end{equation}

\medskip
	\textbf{Claim 3.}~For $\{i,j\}\in {[S]\cup \{0\}\choose 2}$ with $i\neq 0$,  $X_i $ dominates $X_j$ with $ 2(k+1)(r-1) $-deficiency.
\medskip

Fix $i$ and $j$ such that $i\neq 0$.  For every $x\in X_i$, by Claim 1,
\begin{align*}
d_{X_j}(x)&\geq d_{G}(x_i)-(\bigcup_{l\in ([s]\setminus \{i,j\})\cup \{0\}}|X_l|)-|V(\cup_{j=1}^kD_j\cup K)|\\
&=  d_{G}(x_i)-(n-|X_i|-|X_j|)\quad \mbox{(by (\ref{Xi-lower2}))}\\
&\geq |X_j|-(k+1)-(k+1)(2r-3)\\
&=|X_j|-2(k+1)(r-1).
\end{align*}
This completes the proof of Claim 3.

%For every $ x_i\in X_i $ with $ i\ne 0 $, we have
%	\begin{align*}
%		|\overline{N_{G-X_i}(x_i)}|&=(n-d_G(x_i))-|X_i|\\
%		&\leq \left( \frac{1}{r-1}n+(k+1)\right) -\frac{1}{r-1}n+(k+1)(2r-3)\\
%		&=2(k+1)(r-1).
%	\end{align*}
%	Thus $ d_{G-X_i}(x_i)\geq |G-X_i|-2(k+1)(r-1) $ for each $ x_i\in X_i $ where $i=1,...,r-2$. In particular, we have that
%	\begin{equation}\label{2(k+1)(r-1)-defi}
%		d_{X_j}(x)\geq |X_j|-2(k+1)(r-1)
%	\end{equation}
%	for each $ x\in X_i$, i.e., $ X_i $ dominates $X_j$ with $ 2(k+1)(r-1) $-deficiency, where $i=1,...,r-2$, $j=0,1,...,r-2$ and $j\ne i$.
	
	\medskip
	\textbf{Claim 4.}~Let $\{x_1,x_2,...,x_{r-2}\}\subseteq \cup_{i=1}^{r-2} X_i$ such that $x_i\in X_i$ for  $i\in [r-2]$. For any $Y_0\subseteq X_0$ with $|Y_0|\geq 2(k+1)(r-1)^2$,  $ |\cap_{i=1}^{r-2}N_G(x_i)\cap Y_0|\geq k+1 $.
	\medskip
	
	By (\ref{Xi-lower2}),  we may infer that
	$$|\cap_{i=1}^{r-2}N_G(x_i)\cap X_0|\geq |X_0|-2(k+1)(r-1)(r-2).$$
Hence
\begin{align*}
|\cap_{i=1}^{r-2}N_G(x_i)\cap Y_0|&\geq |Y_0|-(|X_0|-|\cap_{i=1}^{r-2}N_G(x_i)\cap X_0|)\\
&\geq 2(k+1)(r-1)^2-2(k+1)(r-1)(r-2)\\
&=2(k+1)(r-1)>(k+1).
\end{align*}
This completes the proof of Claim 4.

	\medskip
	Let $ X_0^* $ denote the set of all vertices of $X_0$ of degree at least $2(k+1)(r-1)^2$ in $X_0$, i.e.,
\[
X_0^*:=\{x\in X_0\ |\ d_{G[X_0]}(x)\geq 2(k+1)(r-1)^2\}.
\]
	
	\medskip
	\textbf{Claim 5.}~$ |X_0^*|\leq 2(k+1)(r-1)(r-2). $
	\medskip
	
	By contradiction, suppose that $ |X_0^*|> 2(k+1)(r-1)(r-2) $. For each $i$, let
\begin{align}\label{X0i-lower}
 X_{0i}^*=\{x\in X_0^*\ |\ d_{X_i}(x)\geq |X_i|/(2(k+1)(r-1)+1)\}.
 \end{align}
	By (\ref{Xi-lower2}), $ d_{X_0}(x_i)\geq |X_0|-2(k+1)(r-1) $ for every $ x_i\in X_i $. Hence for any $ S\subseteq X_0$ with $ |S|=2(k+1)(r-1)+1$, every vertex in $ X_i $ has at least one neighbor in $ S $, which implies that
 \begin{align}\label{Xi-N(S)}
 X_i\subseteq N(S).
 \end{align} By taking the average value, $X_0$ contains at least $|X_0|-2(k+1)(r-1)$ vertices $x$ such that $d_{G[X_i]}(x)\geq 2(k+1)(r-1)^2$. It follows that
 \begin{align}\label{X0i*-lower}
  |X_{0i}^*|\geq |X_0^*|-2(k+1)(r-1).
 \end{align}
By (\ref{X0i*-lower}), we may infer that
\begin{align}\label{X0i-inter-lower}
	 |\cap_{i=1}^{r-2}X_{0i}^*|\geq|X_0^*|-2(k+1)(r-1)(r-2)\geq 1.
\end{align}
   Let $ \mathcal{X}_{0i}^*=X_0^*-X_{0i}^*$.
  For $ T\in {\mathcal{X}_{0i}^*\choose 2(k+1)(r-1)+1} $,  we have $ \sum_{v\in T}d_{X_i}(v)\geq|X_i|$ by (\ref{Xi-N(S)}), which means there exists $v\in T$ such that $d_{X_i}(v)\geq |X_i|/(2(k+1)(r-1)+1)$  contradicting the definition of $\mathcal{X}_{0i}^*$ by (\ref{X0i-lower}).
	
By (\ref{X0i-inter-lower}), there is an $ x_0\in X_0^*  $ such that $ |N(x_0)\cap X_i|\geq |X_i|/(2(k+1)(r-1)+1)$ for each $i=1,...,r-2$.   Since $ n\geq 4(k+1)^2r^4  $, by (\ref{Xi-lower2}) we have
	$$ |N_{X_i}(x_0)|\geq|X_i|/(2(k+1)(r-1)+1)\geq 2(k+1)(r-1)(r-2)+(k+1).$$
Applying Lemma \ref{k cycles} with $ Y_0=N_{X_0}(x_0) $,...,$ Y_{r-2}=N_{X_{r-2}}(x_0) $ and $ b=2(k+1)(r-1) $, we obtain $ k+1 $ vertex disjoint cliques of size $ r-1 $ in $ N(x_0) $. Then $ G $ contains a copy of $ F_{k+1,r} $ contradicting  the fact that $ G $ is $ F_{k+1,r} $-free. This completes the proof of Claim 5.
	
	\medskip
	Let $ Z_0:=X_0-X_0^* $ and $ Z_i:=X_i $ for each $i=1,..,r-2$. By Claim 5 and (\ref{Xi-lower2}), we have that
\begin{align*}
|V(G)-\cup_{i=0}^{r-2}Z_i|&=|V(G)-\cup_{i=0}^{r-2}X_i|+|X_0^*|\\
&\leq (k+1)(2r-3)(r-1)+2(k+1)(r-1)(r-2)\\
&<4(k+1)(r-1)^2.
\end{align*}
Moreover,
\begin{align*}
 |Z_0|&= |X_0|-|X_0^*|\\
 &\geq \frac{n}{r-1}-(k+1)(2r-3)-2(k+1)(r-1)(r-2)\\
 &=\frac{n}{r-1}-2(k+1)(r-1)^2.
	\end{align*}
	
	\medskip
	\textbf{Claim 6.}~For every $v\in V-\cup_{i=0}^{r-2}Z_i$, there exists unique  $ j\in \{0,1,\ldots,r-2\}$ such that $ d_{Z_j}(v)<2(k+1)(r-1)^2+(k+1)$.
	\medskip
	
	Firstly, suppose that there is  $v\in V(G)-\cup_{i=0}^{r-2}Z_i$ such that $ d_{Z_j}(v)\geq2(k+1)(r-1)^2+(k+1)$ for  $j=0,...,r-2$. Let $ b:=2(k+1)(r-1) $ and $ m:=r-1 $. Then for  $j=0,\ldots,r-2$, $d_{Z_j}(v)\geq mb+(k+1)  $ and $ d_{Z_j}(z_i)\geq |Z_j|-b $ for all $ z_i\in Z_i $,  where $ i>0 $ and $ i\ne j $. Applying Lemma \ref{k cycles}, there are $ k+1 $ vertex disjoint cliques of size $ r-1 $ in $ N(x_0) $. Then $ G $ contains a copy of $ F_{k+1,r} $ contradicting the fact that $ G $ is $ F_{k+1,r} $-free.
	
	Secondly, to show the uniqueness, suppose there are two distinct $ j_1 $ and $ j_2 $ such that $ d_{Z_{j_i}}(v)<2(k+1)(r-1)^2+(k+1) $ for both $ i=1,2 $. Since $n\geq 6(k+1)^2r^4 $, we have
	\begin{align*}
		d_G(v)\leq& n-|Z_{j_1}\cup Z_{j_2}|+4(k+1)(r-1)^2+2(k+1)\\
		\leq& n-\left[ (\frac{n}{r-1}-2(k+1)(r-1)^2)+(\frac{n}{r-1}-(k+1)(2r-3))\right]+4(k+1)(r-1)^2+2(k+1)\\
		=&\frac{r-2}{r-1}n-\frac{n}{r-1}+6(k+1)(r-1)^2+(k+1)(2r-1)\\
		<&\frac{r-2}{r-1}n-(k+1).
	\end{align*}
	which contradicts with the fact that $ d_G(x_i)>(\frac{r-2}{r-1})n-(k+1)$.
	This completes the proof of Claim 6.
	
	\medskip
	Adding each $v\in V(G)-(\cup_{i=0}^{r-2}Z_i) $ to $ Z_{j(v)} $, we obtain a partition of $ V(G)=V_0\cup\cdots\cup V_{r-2} $. Clearly, for  $ i\in \{0,..,r-2 \}$,
 \begin{align}\label{Vi-lower}
 |V_i|\geq |Z_i|\geq \frac{n}{r-1}-2(k+1)(r-1)^2.
 \end{align}
  For   $ v_i\in V_i $, by Claims 1 and 6,
\begin{align}\label{Max-degree-Vi}
\Delta(G[V_i])\leq\Delta(G[Z_i])+|V(G)-\cup_{i=0}^{r-2}Z_i|\leq 6(k+1)(r-1)^2.
\end{align}
Now  for $\{i,j\}\in {[r-2]\cup \{0\}\choose 2}$, we have
\begin{align*}
  d_{V_j}(v_i)&=d_G(v_i)-d_{G-V_j}(v_i)\\
  &\geq \delta(G)-(n-|V_j|-(|V_i|-\Delta(G[V_i])))\\
  &\geq |V_j|-(\frac{r-2}{r-1})n-(k+1) -n+(\frac{n}{r-1}-2(k+1)(r-1)^2)-6(k+1)(r-1)^2\quad \mbox{(by (\ref{Vi-lower}) and (\ref{Max-degree-Vi}))}\\
  &>   |V_j|-8(k+1)r^2,
\end{align*}
i.e.,
\begin{align}
 d_{V_j}(v_i)\geq |V_j|-8(k+1)r^2.
 \end{align}
	
	We will show that $ V_0,...,V_{r-2} $ be a partition of $V(G)$ satisfying (i), (ii) and (iii). Let $b':=8(k+1)r^2$. Since $ n\geq 4(k+1)^2r^4 $,  for any $ j $, by (\ref{Vi-lower}), we have
	\[ |V_j|\geq \frac{n}{r-1}-2(k+1)(r-1)^2\geq (r-1)b'+2(k+1) .\]
	
	Now we show $d_{V_i}(y)\leq k$ for all $y\in V_i$.
	Suppose for some $ y\in V_i $, $d_{V_i}(y)\geq k+1 $. Let $ y_1,...,y_{k+1}\in N(y)\cap V_i $ in $ V_i $. By Lemma \ref{tcycles}, there are $ k+1 $ cliques $  D_1,..., D_{k+1} $ such that $ y, y_j\in D_j $ and $ |D_j|=r$ for each $j.$ Further, $ D_j\cap D_l=\{y\} $ for any $ \{j,l\}\in {[k+1]\choose 2} $. $\cup_{i=1}^{k+1}D_i$ is a copy of $F_{k+1,r}$, a contradiction.
	
	Next we show that for all $i\in [r-2]\cup \{0\}$,
\[
\sum_{j\in (\{0\}\cup [r-2])\setminus \{i\} }\nu(G[V_j])\leq  k.
  \]
Suppose that there exists $l\in \{0\}\cup [r-2]$, say $l=0$ such that
\[
\sum_{j\in [r-2]}\nu(G[V_j])\geq k+1.
  \]
 Let $ \{y_1z_1,...,y_{k+1}z_{k+1} \}$ be a $ k+1 $-matching of size $k+1$ in $\cup_{i=1}^{r-2}G[V_i]$. Now, since $ n\geq 4(k+1)^2r^4$,
 \begin{align*}
   |\cap_{i=1}^{k+1}(N_{V_i}(y_j)\cap N_{V_i}(z_j))|&>|V_i|-2(k+1)8(k+1)r^2\\
   &\geq \frac{n}{r-1}-2(k+1)(r-1)^2-16(k+1)^2r^2\geq 1.
 \end{align*}
Therefore, there exists a vertex $ y\in V_0 $ such that $\cup_{j=1}^{k+1}\{y_j,z_j\}\subseteq N(y) $. By Lemma \ref{tcycles}, there are $ k+1 $ cliques $  D_1,..., D_{k+1} $ such that $ y, y_j,z_j\in D_j $ and $ |D_j|=r$ for  $j\in [k+1]$. Furthermore, $ V(D_j)\cap V(D_i)=\{y\} $ for all $ \{j, i\}\in {k+1\choose 2} $. It follows that $G$ contains a copy of $ F_{k+1,r} $, a contradiction.
	
Finally, we show	that (iii) holds by contradiction.
Suppose that there exist $i,v$, say $i=k+1$  such that  $v\in V_{k+1} $  and
\begin{equation}
			d_{G[V_{k+1}]}(x)+\sum_{j=1}^k\nu(G[N(x)\cap V_j])\geq k+1.
		\end{equation}
Write $h_1:=d_{G[V_{k+1}]}(x)$ and $h_2:=\sum_{j=1}^k\nu(G[N(x)\cap V_j])$. Let $v_1,\ldots,v_{h_1}\in N_{G[V_{k+1]}}(x)$ and let
 $ \{y_1z_1,...,y_{h_2}z_{h_2} \}$ be a matching of size $h_2$ in $(\cup_{i=1}^kG[N(x)\cap V_j])$. By (\ref{(1)}), both $ h_1 $ and $ h_2 $ are less than $ k+1 $.  Consider $ k+1 $ of the cliques $ \{v,x_1\},...,\{v,x_s\},\{v,y_1,z_1\},...,\{v,y_t,z_t\} $. Applying Lemma \ref{tcycles} again, there are $ k+1 $ cliques $  D_1,..., D_{k+1} $ in $G$ such that $\{v,x_i\}\subseteq D_i$ for $1\leq i\leq h_1$ and $\{v,y_j,z_j\}\subseteq D_{j+h_1}$ for $1\leq j\leq h_2$. Hence $G$ contains  $ F_{k+1,r} $ as a subgraph, a contradiction. This completes the proof of Lemma \ref{partition}. \qed

	\section{Proof of Theorem \ref{main}}
	\noindent{\textbf{Proof of Theorem \ref{main}.}} By (\ref{low-upp}), we have $ ar(n,F_{k+1,r})\geq\ex(n,F_{k,r})+2 $. So the lower bound is followed. Write $ c(n,k)=\ex(n,F_{k,r})+2. $
	
	For the upper bound, we prove it by contradiction. Suppose that the result does not hold. Then there exists an edge-coloring $ c: E(K_n)\rightarrow [c(n,k)] $ such that $c$ is surjective and the colored graph denoted by $H$ contains no rainbow $ F_{k+1,r} $. Let $G$  be a rainbow subgraph with exactly $c(n,k)$ colors.  Obviously, $G$ is $F_{k+1,r}$-free and $\Delta(G)>(\frac{r-2}{r-1})n$. We discuss two cases. For an edge subset $ Q\subseteq E(H) $, let $ c(Q):=\{c(e)\ |\ e\in Q\} $.
	
	\medskip
	\textbf{Case 1.}~$\delta(G) > (\frac{r-2}{r-1})n-(k+1)$.
	\medskip
	
	By Lemma \ref{partition}, there exists a partition $ V(G)=V_0\cup V_1\cup\cdots\cup V_{r-2} $, such that  (i), (ii) and (iii) hold. By Lemma \ref{edgelow}, we may get that
	\begin{align*}
		\sum_{0\leq i<j\leq r-2}|V_i||V_j|&\geq e(G)-f(k,k)\\
		&=\ex(n,F_{k,r})+2-f(k,k)\\
		&=\begin{cases}
			\ex(n,K_r)-\dfrac{3}{2}k+\dfrac{5}{2},& \mbox{if $k$ is odd},\medskip\\
			\ex(n,K_r)-\dfrac{5}{2}k+2,& \mbox{if $k$ is even}.
		\end{cases}
	\end{align*}
	
	i.e.,
	\begin{align}\label{ViVj-low}
		\sum_{0\leq i<j\leq r-2}|V_i||V_j|&\geq \begin{cases}
			\ex(n,K_r)-\dfrac{3}{2}k+\dfrac{5}{2},& \mbox{if $k$ is odd},\medskip\\
			\ex(n,K_r)-\dfrac{5}{2}k+2,& \mbox{if $k$ is even}.
		\end{cases}
	\end{align}
	
	Let $ t=\max\{||V_j|-\frac{n}{r-1}|, j\in\{0,1,...,r-2\}\} $. Without loss of generality, we may assume $ t=||V_0|-\frac{n}{r-1}| $. Then
	\begin{align*}
		\sum_{0\leq i<j\leq r-2}|V_i||V_j|&\leq |V_0|(n-V_0)+\sum_{1\leq i<j\leq r-2}|V_i||V_j|\\
		&=|V_0|(n-V_0)+\frac{1}{2}\left( (\sum_{j=1}^{r-2}|V_j|)^2-\sum_{j=1}^{r-2}|V_j|^2\right)\\
		&\leq |V_0|(n-V_0)+\frac{1}{2}(n-|V_0|)^2-\frac{1}{2(r-2)}(n-|V_0|)^2\\
		&<-\frac{r-1}{2(r-2)}t^2+\frac{r-2}{2(r-1)}n^2.
	\end{align*}
	where the last second inequality holds by H{\"o}lder's inquality, and the last inequality holds since $ ||V_0|-\frac{n}{r-1}|=t $. That is,
\begin{align}\label{ViVj-upper}
		\sum_{0\leq i<j\leq r-2}|V_i||V_j|<-\frac{r-1}{2(r-2)}t^2+\frac{r-2}{2(r-1)}n^2.
	\end{align}
	By (\ref{ViVj-low}) and  (\ref{ViVj-upper}), we obtain that \begin{align*}
		t^2&<
		\begin{cases}
			\frac{2(r-2)}{r-1}\left( \frac{r-2}{2(r-1)}n^2-\ex(n,K_r)+\dfrac{3}{2}k-\dfrac{5}{2}\right),& \mbox{if $k$ is odd},\medskip\\
			\frac{2(r-2)}{r-1}\left(\frac{r-2}{2(r-1)}n^2-\ex(n,K_r)+\dfrac{5}{2}k-2\right),& \mbox{if $k$ is even}.
		\end{cases}\\
		&<\begin{cases}
			2\left( \frac{r-2}{2(r-1)}n^2-\ex(n,K_r)+\dfrac{3}{2}k-\dfrac{5}{2}\right) ,& \mbox{if $k$ is odd},\medskip\\
			2\left( \frac{r-2}{2(r-1)}n^2-\ex(n,K_r)+\dfrac{5}{2}k-2\right),& \mbox{if $k$ is even}.
		\end{cases}
	\end{align*}
	
	By Tur\'an Theorem, we have $ ex(n,K_r)=\frac{r-2}{2(r-1)}n^2-\epsilon $ where $ \epsilon=\frac{l(r-1-l)}{2(r-1)} $ for $ l\in\{0,1,...,r-2\} $ with $ l\equiv n$   mod $(r-1) $. Note that $ 0\leq \epsilon\leq \frac{r-1}{8} $, and $ \epsilon=0 $ if and only if $ n $ is divisible by $ r-1 $,
	which implies that
	\begin{align*}
		t^2<&\begin{cases}
			2\left( \frac{r-1}{8}+\dfrac{3}{2}k-\dfrac{5}{2}\right) ,& \mbox{if $k$ is odd},\medskip\\
			2\left( \frac{r-1}{8}+\dfrac{5}{2}k-2\right),& \mbox{if $k$ is even}.\end{cases}\\
		=&\begin{cases}
			\frac{r-1}{4}+3k-5, & \mbox{if $k$ is odd},\medskip\\
			\frac{r-1}{4}+5k-4, & \mbox{if $k$ is even}.
		\end{cases}
	\end{align*}
	
	Let $G_i=G[V_i]$ for each $ i=0,...,r-2 $.  Let $G_{cr}$ denote the $ (r-1) $-partite graph with partitions $(V_0,...,V_{r-2})$ and edge set $\cup_{\{i,j\}\in {[r-2]\cup \{0\}\choose 2}}E_G(V_i,V_j)$. By Lemma \ref{partition}, we have
\begin{align}\label{Gi-up-bound}
e(G_i)\leq f(k,k) \quad \mbox{for $i\in \{0,...,r-2\}$.}
\end{align}
    So we may infer that
	\begin{align*}\label{e(G_Cr)-low}
		e(G_{cr})&\geq e(G)-(r-1)f(k,k)\\
		&\geq\begin{cases}
			\ex(n,K_r)+k^2-k+2-(r-1)(k^2+(k-1)/2),&\mbox{if $k$ is odd},\medskip\\
			\ex(n,K_r)+k^2-\frac{3}{2}k+2-(r-1)(k^2+k),& \mbox{if $k$ is even}.
		\end{cases}\\
		&=\begin{cases}
			\ex(n,K_r)-(r-2)k^2-\frac{1}{2}k(r+1)+\frac{1}{2}(r+3),&\mbox{if $k$ is odd},\medskip\\
			\ex(n,K_r)-(r-2)k^2-\frac{1}{2}k(2r+1)+2,& \mbox{if $k$ is even}.
		\end{cases}
	\end{align*}
	
	For $i\in \{0,...,r-2\}$, let $V_{i}':=\{v\in V_i\ |\ N_{G}(v)=V(G)-V_i\}$. Then for $i\in \{0,...,r-2\}$, one can see that the following statements hold:
	\begin{itemize}
		\item [$(i)$]  $|V_i'|\geq |V_i|-(r-2)k^2-\frac{1}{2}k(2r+1)+2$;
		\item [$(ii)$] for any $v\in V_i$, $d_{G_{cr}}(v)\geq n-|V_i|-(r-2)k^2-\frac{1}{2}k(2r+1)+2$.
	\end{itemize}
	
\medskip
\textbf{Claim 5.}~$H[V_i]$ contains no rainbow $K_{1,k+1}$ or a rainbow $(k+1)K_2$.
\medskip

By contradiction. Suppose that $H[V_i]$ contains a rainbow $K_{1,k+1}$ or a rainbow $(k+1)K_2$. Without loss generality, we may assume that $G[V_0]$ contains  $K_{1,k+1}$ or $(k+1)K_2$. Let $R\in \{K_{1,k+1}, (k+1)K_2\}$ be a subgraph of $G$.
% Let $ A=\{e\in E(G)|c(e)\in c(E(K_{1,k+1}))\ \text{or}\ c(E((k+1)K_2))\} $ and $ B=\{e\in E(H)|e\in K_{1,k+1}\ \text{or}\ (k+1)K_2 \} $. As $ G $ is a rainbow graph, we have $ |A|= k+1 $. Define $ V(\hat G)=V(G) $ and $ E(\hat G)=(E(G)\backslash A)\cup B $. 	We suppose $ V(\hat G) $ has the same partitions as $ V(G) $.
%One can see that $ E(\hat G_j)=E(G_j)\backslash A $ for $ j\ne i $ and  $ E(\hat G_i)=(E(G_i)\backslash A)\cup B $, and $ E(\hat G_{cr})=E(G_{cr})\backslash A $.
For any $v\in V_j$ with $j>0$, by $ (ii) $,
\[
d_{G_j}(v_i)\geq |V_j|-(r-2)k^2-\frac{1}{2}k(2r+1)+2-(k+1) \geq |V_j|-2rk^2.
\]
Let $ b:=2rk^2 $. So $ \{V_0,...,V_{r-2}\} $ is $ b $-deficiency complete in $ G$. 	
	Since $ n\geq 16(k+1)^3r^8 $, then  for any $ j $, $ |V_j|> (r-1)b+2(k+1) $.  So by Lemma \ref{tcycles}, $R$ can be greedily extended into a rainbow $F_{k+1,r}$ in $H$, a contradiction. This completes the proof of Claim 5.
	
	By Claim 5, we may assume that for $ i\in\{0,...,r-2\} $, $ H[V_i] $ contains no rainbow $K_{1,k+1}$ or rainbow $(k+1)K_2$. Since $ e(G)>t_{r-1}(n) $, there exists  $i\in \{0,1,...,r-2\}$ such that  $ \nu(G_i)>0 $. 
%\begin{align}\label{Gi-up-bound}
%e(G_i)\leq f(k,k).
%\end{align}
	For $i\in \{0,1,...,r-2\}$, define $S_i=\{ x \in V_i\ |\ N_G(x)=V(G)-V_i \}$ and $s_i=|S_i|$. %That is, $S_i$ is a set of all vertices in $V_i$ which are not incident with every edge in $G_i$ and are adjacent to all vertices in $V(G)-V_{i}$.
 Then, by (\ref{Gi-up-bound}), we have
	\begin{align*}
		\ |S_i|  &\geq |V'_i| - 2(k^2-k+1) \\
&\geq |V_i|-(r-2)k^2-\frac{1}{2}k(2r+1)+2-2(k^2-k+1)\\
&\geq |V_i|\geq \frac{n}{r-1}-(k+1)(2r-3)-rk^2-\frac{1}{2}k(2r+1)+2k\\
		&> 10k^3r^7 \quad\mbox{(since $n\geq 16(k+1)^3r^8$)}.
	\end{align*}
	For $i\in \{0,1,...,r-2\}$, since $H[S_i]$ contains no rainbow matching of size $k+1$ and
	\begin{align*}
		\frac{|S_i|}{2k}  \geq  \frac{ 10k^3r^7}{2k}>5k^2r^7.
	\end{align*}
	$H[S_i]$ has a monochromatic matching $\mathcal{M}_i$ on color $c_i$ with at least $5k^2r^7$ edges.

Let $G'$ be the rainbow subgraph obtained by deleting the  edges  colored by $c_0$ and $c_1$ from $G$. Note that
	\begin{align}
		e(G')\geq \ex(n,F_{k+1,r}).
	\end{align}

\medskip
\textbf{Claim 6.}~$G'$ contains no $F_{k,r}$.
\medskip
	
	Otherwise, suppose that $G'$ contains a copy $R$ of $F_{k,r}$ with center vertex $u$ in $V_z$. We may choose  $ j_0\in\{0,1\} $ such that and $ j_0\ne z $.  Since $\mathcal{M}_{j_0}$ is a monochromatic matching of size at least $5k^2r^7$ and $V(\mathcal{M}_{j_0})\subseteq S_{j_0}$. By the definition of $S_{j_0}$,  there exists an edge $v_{j_0}w_{j_0}\in \mathcal{M}_{j_0}$ such that $uv_{j_0},uw_{j_0}\in E(G')$ and $ \{v_{j_0},w_{j_0}\}\cap V(R)=\emptyset $. By the definition of $G'$, $c(v_{j_0}w_{j_0})\notin \{c(e)\ |\ e\in E(G')\}$. By Lemma \ref{tcycles}, we can greedily find a complete graph $ K$ in $ G'+v_{j_0}w_{j_0}  $ with size of $ r $ such that $V(K)\cap V(R)=\{u\}$. It follows that $K\cup R$ is a rainbow copy of $F_{k+1,R}$ in $H$, a contradiction. This completes the proof of Claim 6.
	
Note that $e(G')\geq ex(n,F_{k,r})$.	By Claim 6,  we have $e(G') = ex(n,F_{k,r})$ and $G'\in EX(n,F_{k,r})$, which also implies that $c_0\neq c_1$. By Theorem \ref{F_{k,r}extremal}, $G_{cr}=G'_{cr}=T(n,r-1)$ and there exists one part set $ V_z $ such that $G'[V_z]\in EX(|V_z|,\{K_{1,k},kK_2\})$ and  $e(G'[V_i]) =0$ for $ i\ne z $.  One can see that $\bigtriangleup(G'[V_z])\leq k-1$. By Theorem \ref{Chvatal},  $G'[V_z]$ contains a matching $M_z$ of size $k-1$.
Recall $|\mathcal{M}_i|\geq 5k^2r^7$ for $i\in \{0,1\}$.  Then there exist $e_0\in \mathcal{M}_0$ and $e_1\in \mathcal{M}_1$ such that $\{e_0,e_1\}\cup M_z$ is a rainbow matching of size $k+1$ in $H$. Furthermore, we may assume  that $e_0,e_1\in E(G)$. Then $\{e_0,e_1\}\cup M_z$ is a  matching of size $k+1$ in $G$. By combining $G_{cr}=T(n,r-1)$, $\{e_0,e_1\}\cup M_z$ can greedily be extended into a copy of $F_{k+1,r}$ in $G$, a contradiction.

 % contradicting Claim 5. %(*****revise CLaim 5)
%Hence, without loss of generality, we may assume that $ z= 0 $ and $ M_z=\{e^1,...,e^{k-1}\} $. Then there exists an edge $e_0\in \mathcal{M}_0$ such that $\{e_0,e^1,...,e^{k-1}\}$ is a rainbow matching of size $k+1$ in $\cup_{i=0}^{r-2}G_i$
	
	%%Recall that $|\mathcal{M}_j|>5k^2r^7$ for $j\in \{0, 1\}$, so we may choose $e_1\in \mathcal{M}_1$ such that  $ M_z\cup\{e_1\} $ is a matching of size $ k $. Since $ V(\mathcal{M}_0)\subseteq S_0 $, by the definition of $ S_1 $, we may choose $ e_0=u_0v_0\in\mathcal{M}_0 $ such that $ u_0 $ and $ v_0 $ are adjacent with all $ V(G)-V_0 $ in $ G^{\prime} $. Then let $ G^{\prime\prime} $ be a graph with vertex set $ V(G^{\prime}) $ and edge set $ E(G^{\prime})\cup \{e_0,e_1\} $.	Note that $ G_{cr}^{\prime\prime}=T(n,r-1) $ and $ G^{\prime\prime}$ is still a rainbow graph. Moreover, we have $ d_{G^{\prime\prime}[V_0]}(u_0)+\sum_{j\ne i}\nu(G^{\prime\prime}[N_{G^{\prime\prime}}(u_0)\cap V_j])=k+1 $. Thus $ G^{\prime\prime} $ contains a rainbow $ F_{k+1,r} $, a contradiction.
	
	\medskip
	\textbf{Case 2.}~$\delta(G) \leq (\frac{r-2}{r-1})n-(k+1)$.
	\medskip
	
	%Suppose that $ n\geq (16(k+1)^3r^8 )^2/(k+1) $. $ G $ is $ F_{k+1,r} $-free with $ c(n,k)=ex(n,K_r)+f(k-1,k-1)+2 $ edges.
Let $ x_0\in V(G) $ such that $ d_G(x_0)=\delta(G)\leq (\frac{r-2}{r-1})n-(k+1) $. Let $ G^0=G $ and $ G^1=G^0-x_0 $. Then by Lemma \ref{delta},
	\begin{equation*}
		e(G^1)=e(G^0-x_0)\geq ex(n-1,K_r)+f(k-1,k-1)+2+(k+1).
	\end{equation*}
	If there exists a vertex $ x_1\in G^1 $ with degree $ d_G(x_1)=\delta(G^1)\leq (\frac{r-2}{r-1})(n-1)-(k+1) $, then delete it to obtain $ G^2=G^1-x_1 $. Continue this process as long as $ \delta(G^t)\leq (\frac{r-2}{r-1})(n-t)-(k+1) $, and after $ l $ steps, we get we get a subgraph $ G^l  $ with $ \delta(G^l)\geq (\frac{r-2}{r-1})l-(k+1) $. By induction, one can see that
	\begin{align*}%\label{l-up-bound}
		ex(n-l,F_{k+1,r})\geq e(G_l)\geq ex(n-l,K_r)+(k+1)l+f(k-1,k-1)+2,
	\end{align*}
which implies that
\begin{align}\label{l-up-bound}
		f(k,k)\geq (k+1)l+f(k-1,k-1)+2.
	\end{align}
By (\ref{l-up-bound}) and Theorem \ref{kk-value}, we may get that $l\leq k$. So we have $ n-l>\sqrt{(k+1)n}\geq 16(k+1)^3r^8  $. Then with the same discussion as Case 1, $ H[V[G^l]] $ contains a rainbow $ F_{k+1,r} $, a contradiction. This completes the proof of Theorem \ref{main}.
	\qed

\end{document}